\documentclass[12pt]{amsart}
\usepackage{amsmath, amssymb, amsthm}
\usepackage{paralist, xcolor, tikz}
\usepackage[noadjust]{cite}
\usepackage{young}
\usepackage[margin=1in,letterpaper,portrait]{geometry}
\usepackage{enumitem}
\usepackage{caption}
\usepackage{comment}
\usepackage{tabularx}
\usepackage{tikz-cd} 
\usepackage[normalem]{ulem}
\usepackage{staves}
\usepackage{hyperref}

\newcommand{\GAV}{\op{GAV}}

\newcommand{\ol}{\overline}

\newcommand{\op}{\operatorname}

\newcommand{\embed}{\iota}
\newcommand{\rc}[1]{#1^{\mathrm{rc}}}

\newcommand{\gl}{\op{gl}}
\newcommand{\mirror}{\vert}
\newcommand{\bigmirror}{\,\big\mirror\,}
\renewcommand{\sb}{S^{\mathrm{B}}}

\newcommand{\note}[1]{{\bf\color{red} #1}}

\newtheorem{conjecture}{Conjecture}[section]
\newtheorem{theorem}[conjecture]{Theorem}
\newtheorem{corollary}[conjecture]{Corollary}
\newtheorem{lemma}[conjecture]{Lemma}
\newtheorem{proposition}[conjecture]{Proposition}

\theoremstyle{definition}
\newtheorem{question}[conjecture]{Question}
\newtheorem{example}[conjecture]{Example}

\newtheorem{definition}[conjecture]{Definition}
\newtheorem*{question*}{Question}

\title{Global patterns in signed permutations}

\author[Owen John Levens]{Owen John Levens$^*$}
\address{Department of Mathematical Sciences, DePaul University, Chicago, IL, USA}
\email{olevens@depaul.edu}
\thanks{$^*$Research partially supported by the DePaul University Undergraduate Summer Research Program.}

\author[Joel Brewster Lewis]{Joel Brewster Lewis$^\dagger$}
\address{Department of Mathematics, George Washington University, Washington, DC, USA}
\email{jblewis@gwu.edu}
\thanks{$^\dagger$Research partially supported by a grant from the Simons Foundation (634530).}

\author[Bridget Eileen Tenner]{Bridget Eileen Tenner$^\ddagger$}
\address{Department of Mathematical Sciences, DePaul University, Chicago, IL, USA}
\email{bridget@math.depaul.edu}
\thanks{$^\ddagger$Research partially supported by NSF Grant DMS-2054436, by a DePaul University Faculty Summer Research Grant, and  by a University Research Council Competitive Research Leave from DePaul University.}

\begin{document}

\begin{abstract}
    Global permutation patterns have recently been shown to characterize important properties of a Coxeter group. Here we study global patterns in the context of signed permutations, with both characterizing and enumerative results. Surprisingly, many properties of signed permutations may be characterized by avoidance of the same set of patterns as the corresponding properties in the symmetric group. We also extend previous enumerative work of Egge, and our work has connections to  the Garfinkle--Barbasch--Vogan correspondence, the Erd\H{o}s--Szekeres theorem, and well-known integer sequences. 
\end{abstract}

\maketitle

\section{Introduction}

In the symmetric group $S_n$, the combinatorics of \emph{pattern avoidance} provides a powerful language for characterizing important families of permutations, as well as a rich source of enumerative questions.  Pattern avoidance has been extended to other Coxeter groups in various ways; recent work by the second and third author revealed the power of \emph{global} permutation patterns as a way to characterize important properties in a Coxeter group \cite{lewis tenner global}. In the present work, we explore global patterns in the context of the hyperoctahedral group of signed permutations. Although the term \emph{global pattern} is new, the idea of global patterns in signed permutations has been independently discovered and studied several times, under several names.  Like classical pattern avoidance in $S_n$, global avoidance provides a powerful language for characterizing signed permutation properties (like \emph{vexillarity} and \emph{smoothness}), as well as a source of interesting enumerative questions---see Section~\ref{sec:survey} for a detailed survey of prior work.

Our goal in the present paper is partly expository.  After introducing the basic notations and terminology of permutations and signed permutations in Section~\ref{sec:main topics}, we compare and contrast the three most studied notions of patterns avoidance/containment for signed permutations (\emph{classical avoidance} of signed permutations by signed permutations and \emph{Billey--Postnikov} and \emph{global} avoidance of permutations by signed permutations) in Section~\ref{sec:three kinds of pattern avoidance}, and survey the state of the art concerning global pattern avoidance in Section~\ref{sec:survey}.

In Section~\ref{sec:characterizing}, we explore several important families of signed permutations (the \emph{layered}, \emph{free}, and \emph{separable} signed permutations, and those with \emph{boolean Bruhat order ideals}) and show that each of them can be characterized by global pattern avoidance.  We observe a striking phenomenon: in all these cases, the class of signed permutations with the property is defined by the avoidance of the \emph{same} patterns as the corresponding class of permutations.  This phenomenon, which we call \emph{persistence}, is remarkable both because of its simplicity and because it does \emph{not} seem to have an analogue for other notions of signed pattern avoidance, where typically a property characterized by avoidance of a small number of patterns in $S_n$ requires many more patterns to characterize in type B and C.

In Section~\ref{sec:enumerations}, we extend prior work of Egge and others on enumerative aspects of global avoidance in signed permutations.  We give a new proof of an old result on the number of signed permutations avoiding the pattern $123$ (equivalently, the pattern $321$), enumerate the classes of signed permutations avoiding a pattern of length $3$ and a monotone pattern, and use the Garfinkle--Barbasch--Vogan correspondence (a signed analogue of Robinson--Schensted) to give an analogue for signed permutations of the Erd\H{o}s--Szekeres theorem. We conclude the paper with several directions for further research in Section~\ref{sec:open}.

\section{The symmetric group, pattern avoidance, and signed permutations}\label{sec:main topics}

The \emph{symmetric group} $S_n$ consists of all permutations of $[n] := \{1, \ldots, n\}$.  We represent an element $w$ of the symmetric group in \emph{one-line notation} as a word $w(1) \cdots w(n)$, or as a permutation matrix in Cartesian coordinates (so the $1$ in the matrix's $i$th column from the left occurs in the $w(i)$th row from the bottom).  

\begin{definition}\label{defn:classical unsigned patterns}
    Fix a permutation $p \in S_k$. A permutation $w \in S_n$ \emph{contains} the \emph{pattern} $p$ if there exist indices $1 \le i_1 < \cdots < i_k \le n$ such that $w(i_1)\cdots w(i_k)$ is order-isomorphic to $p$. Otherwise $w$ \emph{avoids} the pattern $p$.
\end{definition}
This definition can also be naturally expressed in terms of permutation matrices: $w$ contains $p$ if and only if one can choose $k$ rows and $k$ columns of the permutation matrix of $w$ whose intersection is the permutation matrix of $p$.

The symmetric groups belong to an important family of groups called \emph{(finite) Coxeter groups}, which are fully classified; the symmetric groups are \emph{type A} in the classification.  The $n$th \emph{finite Coxeter group of type B}, denoted $\sb_n$, consists of \emph{signed permutations of size $n$}, which we define next.

\begin{definition}\label{defn:signed permutations}
A \emph{signed permutation of size $n$} is a bijection $w$ from the set $[\pm n] := \{\pm1, \pm2, \ldots, \pm n\}$ to itself that satisfies $w(-x) = -w(x)$ for all $x \in [n]$.  A signed permutation in $\sb_n$ may be represented by the $n$-letter \emph{window notation}
$$w(1) \cdots w(n),$$
or by the $2n$-letter \emph{complete notation} as a word
$$w(-n) \cdots w(-1) \bigmirror w(1) \cdots w(n) \ \ = \ \ \big({-w(n)}\big) \cdots \big({-w(1)}\big) \bigmirror w(1) \cdots w(n).$$
\end{definition}
It is customary to save space by writing ``$-x$'' as ``$\ol{x}$,'' and we will often use this convention.  A sequence of integers $w(1) \cdots w(n)$ is the window notation of a signed permutation if and only if $|w(1)|\cdots|w(n)|$ is a permutation of $[n]$.  The reader is referred to \cite{bjorner brenti} as the canonical text on Coxeter combinatorics, and to Chapter 8 of that work, in particular, for its discussion of the group of signed permutations. 

There are several important ``symmetries'' of permutations, of which the most central to our story is the \emph{reverse-complement}. 

\begin{definition}\label{defn:reverse-complement}
    Given a permutation $w \in S_n$, the \emph{reverse-complement} of $w$, denoted $\rc{w}$, is the permutation defined by $\rc{w}(i) = n+1 - w(n+1-i)$ for all $i \in [n]$; that is, in one-line notation,
    $$\rc{w} = \big(n+1 - w(n)\big) \cdots \big(n+1 - w(2)\big) \big(n+1 - w(1)\big).$$
    A permutation $w$ is \emph{invariant} under reverse-complement, or \emph{rc-invariant}, if $w = \rc{w}$.
\end{definition}
Geometrically, the reverse-complement map $p \mapsto \rc{p}$ corresponds to rotation of the graph (equivalently, the permutation matrix) of $p$ by $180^\circ$.

There is a natural embedding $\embed : \sb_n \rightarrow S_{2n}$ defined by sending $w \in \sb_n$ to the unique permutation order-isomorphic to the string
$$w(\ol{n})\,w(\ol{n-1})\,\cdots\, w(\ol{1})\,w(1)\,\cdots \,w(n-1)\,w(n).$$
(This map has appeared under many names in the literature; it was denoted $\iota'$ in \cite{anderson fulton} and $\pi \mapsto \pi^s$ in \cite{egge2007}.) The following result \cite[Prop.~2.2]{egge2007} is clear from the definition.

\begin{proposition}\label{prop:bijection between signed permutations and reverse-complement invariant permutations}
    The map $\embed$ is a bijection from $\sb_n$ to the rc-invariant elements of $S_{2n}$.
\end{proposition}

In the next section, we discuss the extension of pattern avoidance to this context.

\section{Three kinds of signed pattern avoidance}\label{sec:three kinds of pattern avoidance}

In this section, we review the three notions of pattern avoidance in signed permutations that are in widest use: classical avoidance (of a signed permutation pattern), Billey--Postnikov avoidance (of permutation or signed permutation patterns), and global avoidance (of permutation patterns).\footnote{There are also various other notions of pattern avoidance in $\sb_n$ that have arisen and proved useful in specialized contexts, but not been the subject of sustained systematic study, such as the notion appearing in \cite[\S7]{reading} in the context of Cambrian lattices of type B.}

\subsection{Classical pattern avoidance}

If one views a signed permutation as its window notation, it is natural to consider a notion of pattern containment that focuses on the window and considers both the sign and the magnitude of the entries being compared; this leads to the notion of \emph{classical pattern avoidance} for $\sb_n$.

\begin{definition}\label{defn:classical patterns}
Fix a signed permutation $p \in \sb_k$. A signed permutation $w \in \sb_n$ \emph{(classically) contains} the pattern $p$ if there exist indices $1 \le i_1 < \cdots < i_k \le n$ such that
\begin{itemize}
    \item $|w(i_1)| \, \cdots \, |w(i_k)|$ is order-isomorphic to $|p(1)| \, \cdots \, |p(k)|$, and
    \item $w(i_j) \cdot p(j) > 0$ for all $j$.
\end{itemize}
Otherwise $w$ \emph{(classically) avoids} the pattern $p$.
\end{definition}

The same notion arises from considering sub-matrices of signed permutation matrices.  For example, the signed permutation $\ol{4}\,5\,3\,2\,\ol{1}$ classically contains the pattern $\ol{3}\,2\,\ol{1}$ but avoids the pattern $1\,2$; see Figure\,\ref{fig:classical}. 

\begin{figure}[htbp]
\[
\ol{4}\,5\,3\,2\,\ol{1} 
\longleftrightarrow
\begin{pmatrix}
&1&&&\\
\color{red}\mathbf{\ol{1}}&&&&\\
&&\color{red}\mathbf{1}&&\\
&&&1&\\
&&&&\color{red}\mathbf{\ol{1}}
\end{pmatrix}
\hspace{1in}
\ol{3}\,2\,\ol{1}
\longleftrightarrow
\begin{pmatrix}
\ol{1}&&\\
&1&\\
&&\ol{1}
\end{pmatrix}
\]    
 \caption{Classical pattern containment corresponds to taking submatrices: $\ol{4}\,5\,3\,2\,\ol{1}$ contains $\ol{3}\,2\,\ol{1}$ as the submatrix in columns $1, 3, 5$ and rows $4, 3, 1$.  (All unindicated entries are $0$.)}
 \label{fig:classical}
\end{figure}

This notion is perhaps the most widespread version of pattern avoidance in $\sb_n$.  It is not possible to comprehensively chart its appearance in the literature, but it has at least arisen as a tool for characterizing signed permutations with attractive properties in the study of Stanley symmetric functions \cite{billey lam}, the structure of the Bruhat order \cite{tenner patt-bru}, Schubert varieties \cite{billey, lambert}, fully commutative elements \cite{stembridge-fc}, and noncrossing partitions \cite{simion}; it has also been studied from a purely enumerative perspective and from the related perspective of Wilf equivalence \cite{egge2006, mansour-west}.

\subsection{Billey--Postnikov pattern avoidance}

In \cite{billey postnikov}, Billey and Postnikov introduced a notion of pattern containment (now called \emph{Billey--Postnikov patterns}, or \emph{BP-patterns}) that makes sense for elements of any Weyl group.  For our purposes, it is not necessary to give the general root-theoretic definition; instead, we focus on the combinatorial characterization, following \cite[\S2]{haenni}, \cite[\S2]{woo}, and \cite[\S8.3]{woo yong}.\footnote{Because BP-avoidance is defined in terms of root systems rather than directly in terms of their reflection groups, there are actually two versions of BP-avoidance associated to the hyperoctahedral group $\sb_n$, one each from the type-B and type-C root systems. However, it transpires that for two signed permutations $p$ and $w$, one has that $w$ BP-contains $p$ as type-B objects if and only if $w$ BP-contains $p$ as type-C objects, and furthermore that no type-B object contains a type-C object or vice-versa. 
Thus the distinction is ultimately without a difference from our perspective.}

In the case of an element $w$ of the hyperoctahedral group $\sb_n$, the patterns it may contain are signed permutations in $\sb_k$ for $k \leq n$ and permutations in $S_k$ for $k \leq n$.  In the case of a pattern that is also of type B, BP-containment coincides precisely with classical pattern containment, as in Definition~\ref{defn:classical patterns}.  In the case of a type-A pattern, BP-containment can be characterized in the following way.

\begin{definition}\label{defn:type-A BP-avoidance}
Fix a permutation $p \in S_k$.  A signed permutation $w \in \sb_n$ \emph{BP-contains} the pattern $p$ if there exist indices $i_1 < \cdots < i_k$ in $[\pm n]$ such that
\begin{itemize}
    \item for $a \in [\pm n]$, if $a \in \{i_1, \ldots, i_k\}$ then $\ol{a} \not \in \{i_1, \ldots, i_k\}$, and
    \item $w(i_1) \cdots w(i_k)$ is order-isomorphic to $p$.
\end{itemize}
Otherwise, $w$ \emph{BP-avoids} the pattern $p$.
\end{definition}

For example, the signed permutation $w = \ol{2}\,1$ BP-contains the pattern $p = 21 \in S_2$ because $w(\ol{2}) = \ol{1} > \ol{2} = w(1)$ and the index set $\left\{\ol{2}, 1\right\}$ does not contain both $a$ and $\ol{a}$ for any $a \in [\pm n]$, but the signed permutation $v = \ol{1}\,2$ BP-avoids the pattern $p$, because the only inversion in the complete notation $\ol{2}\,1 \bigmirror \ol{1}\,2$ for $v$ is between the entries in positions $\ol{1}$ and $1$.  (Both $w$ and $v$ BP-avoid the pattern $21 \in \sb_2$.)

BP-avoidance has proven extremely useful in the characterization of numerous signed permutation properties: these include the elements that index smooth or rationally smooth Schubert varieties \cite{billey postnikov}, the elements having the property that the Bruhat ideal they generate has the same cardinality as the set of regions formed by their inversion arrangement (among other characterizations) \cite{woo}, the fully commutative elements \cite{billey postnikov}, the elements with boolean Bruhat ideals \cite{gao hanni boolean}, and the separable elements \cite{gaetz gao}, among others.

\subsection{Global pattern avoidance}\label{sec:global}

A third notion of pattern avoidance by signed permutations is common in the literature.  It historically does not have a standard name; the present authors hope that the name \emph{global avoidance} will catch on.\footnote{See also \cite{W}.} 

\begin{definition}\label{defn:global patterns}
Fix a permutation $p \in S_k$. A signed permutation $w \in \sb_n$ \emph{globally contains} the pattern $p$ if there exist indices $i_1 < \cdots < i_k$ in $[\pm n]$ such that $w(i_1)\cdots w(i_k)$ is order-isomorphic to $p$. Otherwise $w$ \emph{globally avoids} the pattern $p$. 
\end{definition}

\begin{example}\label{ex:global containment and avoidance}
Consider the signed permutation $w = \ol{2}\,1\,3\,\ol{4} = 4 \,\ol{3} \,\ol{1} \,2 \bigmirror \ol{2}\,1\,3\,\ol{4} \in \sb_4$. 
\begin{enumerate}[label = (\alph*)]
    \item Every element of $S_3$ is globally contained in $w$, as marked in red in Table~\ref{tab:s3 in -2 1 3 -4}.
    \begin{table}[htbp]
    ${\renewcommand{\arraystretch}{1.5}\begin{array}{c|r@{\bigmirror}l}
        \text{Pattern $p$} & \multicolumn{2}{c}{\text{Global occurrence of $p$ in $w$}}\\
        \hline
        123 & \hspace{.5in} 4 \, \note{\ol{3}} \,\ol{1}\, 2\! &         \ol{2}\,\note{1} \, \note{3}\,\ol{4}\\
        132 & 4 \,\note{\ol{3}}\, \ol{1}\, \note{2} & \ol{2}\,\note{1}\,3\,\ol{4}\\
        213 & 4 \,\ol{3} \,\note{\ol{1}}\, 2 & \note{\ol{2}}\,1\,\note{3}\,\ol{4}\\
        231 & 4\, \ol{3} \,\note{\ol{1}}\, \note{2} & \note{\ol{2}}\,1\,3\,\ol{4}\\
        312 & \note{4} \,\ol{3}\, \ol{1}\, 2 & \ol{2}\,\note{1}\,\note{3}\,\ol{4}\\
        321 & \note{4} \,\ol{3} \,\ol{1} \,\note{2} & \ol{2}\,\note{1}\,3\,\ol{4}
    \end{array}}$
    \caption{Global occurrences of elements of $S_3$ in $\ol{2}\,1\,3\,\ol{4}$.}
    \label{tab:s3 in -2 1 3 -4}
    \end{table}
    
    \item The signed permutation $w$ globally avoids the pattern $2143$. Indeed, the only global occurrences of $213$ in $w$ are $\ol{1} \,\ol{2} \,1$, $\ol{1}\, \ol{2} \,3$, $2\,\ol{2}\,3$, and $213$, none of which can be extended to a global $2143$-pattern in $w$.
\end{enumerate}
\end{example}

We denote by $\GAV_n(P)$ the set of elements of $\sb_n$ that globally avoid the set $P$ of patterns, and by $\GAV(P)$  the union of such sets over all $n$. When $P = \{p\}$ is a singleton set, we write $\GAV_n(p)$ and $\GAV(p)$. 

This notion of pattern avoidance is the one that arises automatically from the embedding $\embed: \sb_n \to S_{2n}$.

\begin{proposition}\label{prop:bijection between global avoidance and classical avoidance in reverse-complement invariant permutations}
 A signed permutation $w$ globally contains a permutation $p$ if and only if the permutation $\embed(w)$ contains $p$ (in the sense of Definition~\ref{defn:classical unsigned patterns}).
\end{proposition}

Because the image $\embed(w)$ of every signed permutation $w$ is fixed by reverse-complement, certain redundancies follow in the study of global patterns.
  
For any dihedral symmetry $s$ of a square and any permutation $p$, let $p^s$ denote the permutation whose graph is the result of applying $s$ to the graph of $p$. 
For classical pattern avoidance in $S_n$, it is easy to show that the patterns $p$ and $p^s$ are avoided by equally many permutations of any fixed size; indeed, a permutation $q$ avoids $p$ if and only if $q^s$ avoids $p^s$. This phenomenon extends to $\sb_n$ in a straightforward way.

\begin{lemma}[{\cite[Lem.~2.6]{egge2007}}]
\label{lem:egge symmetries}
Let $P$ be any set of patterns. 
\begin{enumerate}[label = (\alph*)]
\item Let $s$ be a dihedral symmetry and let $P^s$ be the result of applying $s$ to every element of $P$.  Then $|\GAV_n(P)| = |\GAV_n(P^s)|$.
\item Let $Q$ be a set that contains, for each pattern $p \in P$, either $p$ or $\rc{p}$ or both (and no other patterns).  Then $\GAV(P) = \GAV(Q)$.
\end{enumerate}
\end{lemma}

One consequence of Lemma~\ref{lem:egge symmetries}(b) is that, unlike classical pattern classes, most global pattern classes are not defined by a unique minimal set of avoided patterns. For example, $\GAV(132) = \GAV(213)$ and $\GAV(1) = \GAV(\{12, 21\})$.

Global avoidance has been explored both for its use in characterizing various signed permutation properties and from a purely enumerative/combinatorial point of view; we defer a more detailed survey of this prior work to Section~\ref{sec:survey}.

\subsection{Compare and contrast}

We begin by making some simple observations about differences between the three notions of pattern containment in Definitions~\ref{defn:classical patterns}, \ref{defn:type-A BP-avoidance}, and~\ref{defn:global patterns}.  For any version of pattern containment, say that a subset of $\bigcup_{n \geq 0} \sb_n$ is a \emph{pattern class} if it consists precisely of the (signed) permutations that avoid a fixed set of patterns. A first observation is that the patterns in classical signed containment are signed permutations, while those in type-A BP-containment and global containment have unsigned permutations as patterns.  Thus, aspects of pattern containment that relate to the \emph{pattern poset} do not apply to these settings.  One concrete consequence is that although classical pattern classes have unique bases (i.e., if the set of signed permutations that classically avoid all patterns in $P$ is the same as the set of signed permutations that classically avoid all patterns in $Q$, then $P = Q$), the same is not true for either BP-classes or global classes. 
A second observation is that the indices $i_1, \ldots, i_k$ at which the pattern appears in Definitions~\ref{defn:classical patterns} and~\ref{defn:type-A BP-avoidance} cannot include any pair $\{\pm i\}$, while in Definition~\ref{defn:global patterns} this restriction does not hold.  As a result, a signed permutation $w\in \sb_n$ may globally contain permutations of size up to $2n$, whereas in classical containment and BP-containment the size of contained patterns is bounded by $n$. 

A third comparison concerns how well the three notions generalize.  BP-avoidance makes sense for any Weyl group, by definition; however, it has the odd feature that BP-avoidance in type A (of permutations by other permutations) does not exactly coincide with the classical notion of pattern avoidance: for example, the set of $132$-avoiding permutations is not a BP-pattern class for any set of BP-patterns.  Global containment naturally generalizes to any of the \emph{George groups}, the finite and affine Coxeter groups of types A, B, C, D; in the case of type A, it coincides with classical pattern containment.

Another natural way to compare two notions of pattern avoidance is in terms of their ``expressiveness,'' by asking whether every pattern class of one kind is also a pattern class of the second kind.  Since every type-A BP-pattern class and every global pattern class is a classical class \cite[Lem.~2.29]{haenni} \cite[Thm.~2.20]{egge2007}, and since the classical class of signed permutations avoiding $\ol{1}$ is neither a global pattern class nor a type-A BP-pattern class,\footnote{Because, under the natural inclusions $S_n \subset \sb_n$, every permutation is an element of this set, and so every permutation is BP- and globally contained in some element of this set; but the only BP- or global pattern class with this property is the whole hyperoctahedral group $\sb_n$.}
we have that classical avoidance is strictly more expressive than both type-A BP-avoidance and global avoidance. We discuss the relative expressiveness of these latter two notions now. The next example shows that type-A BP-avoidance is \emph{not} capable of capturing global avoidance in general.

\begin{example}
The class of signed permutations globally avoiding $12$ is not a type-A BP-avoidance class.  This is because $\GAV_1(\{12\}) = \left\{\ol{1}\right\}$ contains only one of the two signed permutations in $\sb_1$, but every type-A BP-pattern class intersects $\sb_1$ in the empty set (when BP-avoiding a set containing $p=1$) or $\{1,\ol{1}\}$ (when BP-avoiding anything else). 
\end{example}

This leaves open the question of whether global avoidance is strictly more expressive than type-A BP-avoidance, or whether the two notions are incomparable.  To answer this question, it would suffice to either produce a single permutation $p$ whose BP-avoidance class is not a global avoidance class (of any set of permutations), or to show that the BP-avoidance classes of individual permutations are all global avoidance classes. For every pattern $p$ in $S_1$, $S_2$, $S_3$, and for most patterns in $S_4$, we have been able to show by construction that the BP-class avoiding $p$ is a global avoidance class.  In each case, the global class uses patterns that are exactly twice as long as $p$.  We describe the exception now.

\begin{example}
Let $p = 1 4 3 2 \in S_4$.  Then $p$ is BP-contained in the signed permutation in $1\,\ol{4}\,3\,\ol{2} \in \sb_4$, namely, in positions $-3$, $-2$, $1$, $4$ of its complete notation $2\,\note{\ol{3}}\,\note{4}\,\ol{1} \bigmirror \note{1}\,\ol{4}\,3\,\note{\ol{2}}$.
Let $q = \embed( 1\,\ol{4}\,3\,\ol{2} ) = 6 2 8 4 5 1 7 3$ be the image of this signed permutation in $S_8$.  If the BP-class avoiding $p$ were characterized by the global avoidance of patterns of size $8$, then necessarily $q$ would be one of these patterns (as it is the only permutation in $S_8$ that is globally contained in $1\,\ol{4}\,3\,\ol{2}$), and therefore any signed permutation that globally contained $q$ would necessarily BP-contain $p$.  However, the signed permutation $w = 5\, \ol{3}\, \ol{2}\, \ol{6}\, 1\, \ol{4} \in \sb_6$ globally contains $q$, at positions $-5$, $-1$, $1$, $2$, $3$, $4$, $5$, $6$ in its complete notation
\[
4 \, \note{\ol{1}} \, 6 \, 2 \, 3 \, \note{\ol{5}} \bigmirror \note{5 \, \ol{3} \, \ol{2} \, \ol{6} \, 1 \, \ol{4}},
\]
while one can check that $w$ does \emph{not} BP-contain $p$.
\end{example}

This example naturally raises the following questions, which we have not managed to resolve.

\begin{question}
Is the set of signed permutations that BP-avoid the pattern $1432 \in S_4$ a global avoidance class for any set of patterns?  If so, is \emph{every} BP-avoidance class a global avoidance class?
\end{question}

\section{A survey of global pattern avoidance}
\label{sec:survey}

The remainder of this paper is focused on global pattern avoidance. In this section, we endeavor to give a reasonably comprehensive account of prior appearances in the literature of global pattern containment.   
We begin with several important families of signed permutations that are characterized by global avoidance.

\subsection{Vexillarity}\label{sec:vex}

\emph{Vexillary permutations} were introduced by Lascoux and Sch\"utzenberger \cite{lascoux schutzenberger} in their study of Schubert polynomials, and discovered independently by Stanley in \cite{stanley 1984}.  They have a number of equivalent definitions (see \cite[(1.27)]{macdonald}), which include the following two: a permutation $w \in S_n$ is vexillary if and only if $w$ avoids $2143$, if and only if its Stanley symmetric function is a single Schur function. 
In \cite{billey lam}, Billey and Lam generalized the second notion to signed permutations, defining a signed permutation $w \in \sb_n$ to be vexillary if its Stanley symmetric function is a single Schur $Q$-function.
Billey and Lam then showed \cite[Thm.~14]{billey lam} that these vexillary signed permutations can be characterized by the classical avoidance of the nine signed patterns 
\[
\left\{
    21, \ 
    \ol{3}\,2\,\ol{1}, \ 
    2\,\ol{3}\,4\,\ol{1}, \
    \ol{2}\,\ol{3}\,4\,\ol{1}, \ 
    3\,\ol{4}\,\ol{1}\,\ol{2}, \
    \ol{3}\,\ol{4}\,1\,\ol{2}, \
    \ol{3}\,\ol{4}\,\ol{1}\,\ol{2}, \
    \ol{4}\,1\,\ol{2}\,3, \
    \ol{4}\,\ol{1}\,\ol{2}\,3
    \right\}
\]
(entry P0032 of \cite{dppa}).  
In \cite{anderson fulton}, Anderson and Fulton showed that this property can be neatly described in terms of global pattern avoidance: their main theorem asserts (in part) that a signed permutation is vexillary for type B if and only if it globally avoids the pattern $2143$.

\subsection{Full commutativity}\label{sec:fully commutative}

A permutation $w$ is \emph{fully commutative} if any reduced decomposition for $w$ can be obtained from any other by a series of commutations between consecutive generators (i.e., without using any longer braid moves), or equivalently if no reduced word for $w$ contains consecutive letters $i(i \pm 1)i$.  The fully commutative permutations are precisely the permutations that avoid the pattern $321$ \cite{billey jockusch stanley}.  Several different extensions of this notion to signed permutations have been defined.  In particular, Stembridge showed \cite[Cor.~5.4]{stembridge-fc} that every fully commutative signed permutation belongs to one of two classes he called \emph{top} and \emph{bottom}; the fully commutative top elements are the same as the signed permutations that Fan \cite{fan} called \emph{commutative elements of type C}.  These signed permutations are precisely those that globally avoid $321$. 
See \cite{lewis tenner global} for more detailed discussion of this case, which also extends to the affine symmetric group \cite{Green}.

\subsection{Smoothness}\label{sec:smooth}

For background on this section, see \cite{abe billey} and \cite[Ch.~3]{manivel}.

The (type-A) \emph{Schubert varieties} are indexed by permutations.  A natural problem is to determine which Schubert varieties are (rationally) smooth, in terms of the indexing permutation.  This characterization was carried out by Lakshmibai and Sandhya \cite{lakshmibai sandhya}, who showed that a permutation $w \in S_n$ indexes a smooth Schubert variety of type A if and only if $w$ avoids the patterns $3412$ and $4231$ (entry P0005 of \cite{dppa}).

Schubert varieties can also be defined for other types; in types B and C, they are indexed by signed permutations.  A classification of the $w \in \sb_n$ for which the associated type-C Schubert variety is smooth was carried out in \cite[Thm.~1]{lakshmibai song}, in terms of several criteria that include a global pattern criterion.  In particular, \cite[Thm.~2]{lakshmibai song} asserts that if $w$ globally avoids $3412$ and $4231$, then it indexes a smooth type-C Schubert variety.

The characterization of smooth and rationally smooth  Schubert varieties in types B and C was carried out by Billey in terms of classical pattern avoidance \cite[Thm.~6.1]{billey}.  These criteria were reformulated in terms of BP-avoidance in \cite{billey postnikov}, and in terms of global pattern avoidance (under the name ``extended patterns'') in \cite[\S\S8 and 13.3]{billey lakshmibai}.\footnote{The discussion in \cite{billey lakshmibai} on this point is not fully explicit.  However, in a personal communication, Sara Billey explained that she and Lakshmibai checked (in unpublished work) that the work of \cite{lakshmibai song} can be extended to show that these are global pattern classes. By a computation from the classical pattern characterization, one can then show that, in particular, the global avoidance of the patterns in \cite[\S13.3]{billey lakshmibai} also characterizes this set.} An immediate corollary of this work is the following result.

\begin{corollary}[{\cite[Prop.~3.2]{skandera}}]\label{thm:smooth is cool}
    The set $\GAV_n(\{3412, 4231\})$ of signed permutations in $\sb_n$ that globally avoid the patterns $3412$ and $4231$ is equal to the set of signed permutations that simultaneously index smooth Schubert varieties of both types B and C.
\end{corollary}

\subsection{Enumeration}\label{subsec:enum}

In \cite{egge2007, egge2010}, Egge made a detailed study of pattern avoidance in the classes of permutations fixed by some subgroup of dihedral symmetries.  In particular, Egge's results on rc-invariant permutations translate immediately into results on global pattern avoidance for signed permutations.  In particular, Egge enumerated signed permutations globally avoiding subsets of patterns of length $3$ \cite[Thm.~2.10]{egge2007}, with the most involved case being a recursive/generatingfunctional argument to show that there are $\binom{2n}{n}$ elements of $\sb_n$ that globally avoid $123$ \cite[Thm.~2.17]{egge2007}, and in \cite{egge2010} he used the Garfinkle--Barbasch--Vogan (GBV) correspondence between signed permutations and pairs of domino tableaux \cite{garfinkle, barbasch vogan, van leeuwen} to enumerate signed permutations globally avoiding the monotone pattern $12\cdots k$ for any $k$. 
Subsequent work by Lonoff and Ostroff \cite{lonoff ostroff} and Troyka \cite{troyka} enumerated rc-invariant permutations avoiding pairs of short patterns, and Troyka explored the relationship between the growth rates of a permutation class and the corresponding signed permutation global class.

Other enumerative work has involved the global classes mentioned in the preceding subsections.  Gao and H\"anni \cite{gao hanni vex} showed that $|\GAV_n(2143)| = |\GAV_n(1234)|$ for all $n$, answering a question of Anderson and Fulton \cite{anderson fulton}, and the $321$-avoiding signed permutations are enumerated by the central binomial coefficient $\binom{2n}{n}$---see Section~\ref{subsec:321} for more on this case. 

The smooth and rationally smooth Schubert varieties of types B and C are enumerated in \cite{richmond slofstra}.  In light of the discussion in Section~\ref{sec:smooth}, these results may be viewed as enumerations of global pattern classes.  In particular, the authors obtain the surprising result that the number of signed permutations in $\sb_n$ indexing smooth Schubert varieties in both types B and C simultaneously is the same as the number of permutations in $S_{n + 1}$ indexing smooth Schubert varieties of type A, i.e., that the number of globally $3412$- and $4231$-avoiding signed permutations in $\sb_n$ is equal to the number of $3412$- and $4231$-avoiding permutations in $S_{n + 1}$.

\section{Characterization results}\label{sec:characterizing}

As we have already seen, an important and influential portion of combinatorial research in permutation patterns concerns properties that can be characterized by pattern avoidance (as catalogued in \cite{dppa}). In many cases, a property that can be characterized in terms of (classical) pattern avoidance in $S_n$ (i.e., the permutations that possess the property form a pattern class) has been generalized to signed permutations in $\sb_n$, with its own characterization by classical signed patterns. Put another way, one has the situation that
\begin{align*}
    \{w \in S_n : w \text{ has property } \phi\} &= \{w \in S_n : w \text{ avoids all patterns in } A_{\phi}\} \ \  \text{ and}\\
    \{w \in \sb_n : w \text{ has property } \phi\} &= \{w \in \sb_n : w \text{ classically avoids all patterns in } B_{\phi}\}
\end{align*}
for some sets $A_\phi$ and $B_\phi$ of patterns. Typically $|B_{\phi}| > |A_{\phi}|$; for example, we see this in the case of vexillary permutations \cite[P0004]{dppa} described above in Section~\ref{sec:vex}.  A similar phenomenon is common in the case of BP-avoidance: signed permutation properties seem to usually have more complicated characterizations as BP-classes than the corresponding (unsigned) permutation classes do.

In this section, we show two surprising phenomena: first, that many important pattern classes of signed permutations can be expressed in terms of global (not just classical) pattern avoidance, and second, that when we do this, it often transpires that the \emph{same} global patterns characterize the type-A and type-B classes.  Therefore we make the following definition.

\begin{definition}
    Let $\phi$ be a property defined for elements of $S_n$ and of $\sb_n$, such that 
    $$\{w \in S_n : w \text{ has property } \phi\} = \{w \in S_n : w \text{ avoids all patterns in } A_{\phi}\}$$
    for a set $A_{\phi}$ of unsigned permutations. If 
    $$\{w \in \sb_n : w \text{ has property } \phi\} = \{w \in \sb_n : w \text{ globally avoids all patterns in } A_{\phi}\},$$
    then the property $\phi$ is \emph{persistent}. 
\end{definition}

In this language, it has already been observed that the vexillary property is persistent, while some related phenomena were observed for the properties of being smooth or fully commutative.  Below, we consider several other interesting cases, and show that the permutation properties \emph{boolean}, \emph{free}, \emph{separable}, and \emph{layered} are all persistent to type B.  On the other hand, being \emph{Grassmannian} or \emph{bigrassmannian}, or having \emph{rank-symmetric} or \emph{self-dual} Bruhat intervals, are not persistent.

\subsection{Boolean Bruhat order ideals}

In this and the next section, we discuss the \emph{Bruhat order} of a Coxeter group $G$ and the \emph{reduced decompositions} of the elements of $G$.  The definitive reference for these topics is \cite[Chs.~1--2]{bjorner brenti}.
The Bruhat order is a  notoriously complex and interesting partial ordering, and combinatorial aspects of the poset translate to important algebraic and geometric properties of the Coxeter group. In this ordering on a given group, there are some elements whose principal order ideals are particularly nice: the \emph{boolean} elements.

\begin{definition}\label{defn:boolean}
    Given a Coxeter group $G$, an element $w \in G$ is \emph{boolean} if the principal order ideal of $w$ in the Bruhat order on $G$ is isomorphic to a boolean lattice.
\end{definition}

Boolean elements can be characterized by a feature of their reduced decompositions.

\begin{proposition}[{\cite[Prop.~7.3]{tenner patt-bru}}]
    Given a Coxeter group $G$, an element $w \in G$ is boolean if and only if the reduced decompositions of $w$ contain no repeated generators.
\end{proposition}

The boolean elements in the finite Coxeter groups of types A, B, and D were characterized in \cite{tenner patt-bru} in terms of classical pattern avoidance. 

\begin{proposition}[{\cite[Thms.~4.3 and~7.4]{tenner patt-bru}}]
    A permutation in $S_n$ is boolean if and only if it classically avoids the patterns $321$ and $3412$.
    A signed permutation in $\sb_n$ is boolean if and only if it classically avoids the signed patterns 
    \begin{equation}\label{eqn:type b boolean}
        \left\{1\,\ol{2}, \ \ol{1}\,\ol{2}, \ \ol{2}\,\ol{1}, \ 321, \ 3\,2\,\ol{1}, \ 3\,\ol{2}\,1, \ \ol{3}\,2\,1, \ 3412, \ 3\,4\,\ol{1}\,2, \ \ol{3}\,4\,1\,2\right\}.
    \end{equation}
\end{proposition}

This was subsequently reformulated in terms of  BP-avoidance in \cite{gao hanni boolean}: the boolean signed permutations are those that BP-avoid the type-A patterns $321$ and $3412$ and also the type-B/C patterns $1\, \ol{2}$, $\ol{2}\, \ol{1}$, $\ol{1}\, \ol{2}$, and $3 \, \ol{2} \, 1$.  It turns out that this property can be expressed more elegantly in terms of global avoidance---it is persistent.

\begin{theorem}\label{thm:type b boolean global patterns}
    A signed permutation is boolean for type B if and only if it globally avoids the patterns $321$ and $3412$. In other words, being boolean is a persistent property.
\end{theorem}

\begin{proof}
    Suppose, first, that $w$ classically contains one of the ten patterns listed in~\eqref{eqn:type b boolean}. Looking at the complete notation, as marked below, each of these patterns globally contains a $321$- or $3412$-pattern (and often more than one). 
    $$\begin{array}{lr@{\bigmirror}lcr@{\bigmirror}lcr@{\bigmirror}l}
        & 2 \, \note{1} & \note{\ol{1}} \, \note{\ol{2}}
        & \hspace{.25in} &
        \ol{1} \, \ol{2} \, \ol{3} & \note{3} \, \note{2} \, \note{1}
        & \hspace{.25in} &
        \ol{2} \, \ol{1} \, \ol{4} \, \ol{3} & \note{3} \, \note{4} \, \note{1} \, \note{2}\\
        \\
        & \note{1} \,\note{2} & \note{\ol{2}} \,\note{\ol{1}}
        & &
        1 \, \ol{2} \, \ol{3} & \note{3} \, \note{2} \, \note{\ol{1}}
        & &
        \ol{2} \, 1 \, \ol{4} \, \ol{3} & \note{3} \, \note{4} \, \note{\ol{1}} \, \note{2}\\
        \\
        & \note{2} \, \ol{1} & \note{1} \, \note{\ol{2}}
        & &
        \ol{1} \, \ol{2} \, \note{3} & \ol{3} \, \note{2} \, \note{1}
        & &
        \ol{2} \, \ol{1} \, \ol{4} \, \note{3} & \ol{3} \, \note{4} \, \note{1} \, \note{2}\\
        \\
        \multicolumn{4}{c}{ } & \ol{1} \, \note{2} \, \ol{3} & \note{3} \, \note{\ol{2}} \, \note{1}
    \end{array}$$
 Thus $w$ globally contains a $321$- or $3412$-pattern.

    Conversely, consider a signed permutation $w \in \sb_n$ that globally contains $321$ or $3412$. 
    Suppose, first, that $w$ globally contains $321$, and that this global pattern occurrence is in positions $p_1 < p_2 < p_3$, meaning that $w(p_3) < w(p_2) < w(p_1)$. Our goal is to show that this $w$ must classically contain one of the patterns listed in~\eqref{eqn:type b boolean}.

    We will proceed by cases, based on the signs of the positions $p_i$. 
    For brevity, we leave some details as an exercise to the reader, but we outline the results here. 
    \begin{itemize}
        \item When $p_1 > 0$ or $p_3 < 0$, the signed permutation $w$ classically contains at least one element of
        $$\left\{ \ol{1}\,\ol{2}, \ 
        321, \
        3\,2\,\ol{1}, \ 
        1\,\ol{2}
        \right\}.$$
        \item When $p_1 < 0 < p_2$ or $p_2 < 0 < p_3$, the signed permutation $w$ classically contains at least one element of 
        $$\left\{
        \ol{1}\,\ol{2}, \
        \ol{2}\,\ol{1}, \ 
        \ol{3}\,2\,1, \ 
        1\,\ol{2}
        \right\}.$$
    \end{itemize}

    Now suppose that $w$ globally contains $3412$, and that this global pattern occurrence is in positions $p_1 < p_2 < p_3 < p_4$, meaning that $w(p_3) < w(p_4) < w(p_1) < w(p_2)$. We again proceed by cases.
    \begin{itemize}
        \item When $p_1 > 0$ or $p_4 < 0$, the signed permutation $w$ classically contains at least one element of
        $$\left\{\ol{2}\,\ol{1}, \ 
        3412, \ 
        3\,4\,\ol{1}\,2, \
        1\,\ol{2}, \ 
        3\,\ol{2}\,1
        \right\}.$$
        \item When $p_1 < 0 < p_2$ or $p_3 < 0 < p_4$, the signed permutation $w$ classically contains at least one element of
        $$\left\{\ol{1}\,\ol{2}, \
        \ol{2}\,\ol{1}, \ 
        \ol{3}\,4\,1\,2, \ 
        1\,\ol{2}, \ 
        3\,\ol{2}\,1
        \right\}.$$
        \item When $p_2 < 0 < p_3$, the signed permutation $w$ classically contains the pattern $\ol{2}\,\ol{1}$.
    \end{itemize}

   Thus, in whatever way the signed permutation $w$ globally contains a $321$- or $3412$-pattern, it also classically contains at least one of the patterns listed in~\eqref{eqn:type b boolean}.

    We have now shown that global containment of $\{321,3412\}$ is equivalent to classical containment of at least one pattern listed in~\eqref{eqn:type b boolean}. Therefore global avoidance of $321$ and $3412$ is equivalent to classical avoidance of those ten patterns, meaning that global $\{321,3412\}$-avoidance characterizes the boolean elements of type B.
\end{proof}

\subsection{Freeness}

Another persistent pattern property is that of \emph{freeness}, introduced in \cite{petersen tenner depth}.

\begin{definition}\label{defn:free}
    A permutation is \emph{free} if, in its reduced decompositions, all generators commute with each other.
\end{definition}

Freeness can be characterized by pattern avoidance, and appears as \cite[P0026]{dppa}.

\begin{proposition}[{\cite{petersen tenner depth}}]\label{prop:free patterns in type a}
The free permutations in $S_n$ are exactly the permutations avoiding the patterns $\{231,312,321\}$.
\end{proposition}

Definition~\ref{defn:free} makes sense in the setting of signed permutations, and so we have free elements of $\sb_n$, as well. Moreover, this description can be rephrased in terms of classical signed pattern avoidance. We follow \cite[Ch.~8]{bjorner brenti} for notation of the simple generators in $\sb_n$.

\begin{proposition}\label{prop:characterizing free signed perms}
    The free elements of $\sb_n$ are the signed permutations that classically avoid $
    \left\{
        231, \
        312, \
        321, \
        \ol{2}\,1, \
        \ol{1}\,\ol{2}, \
        \ol{2}\,\ol{1}, \
        2\,\ol{1}, \
        1\,\ol{2}
    \right\}$.
\end{proposition}

\begin{proof}
    Consider a free element $w \in \sb_n$, and let $\op{supp}(w) \subseteq [0,n-1]$ be the \emph{support} of $w$. That is, $i \in \op{supp}(w)$ if and only if the generator $\sb_i$ appears in reduced decompositions of $w$. To be free means that $\op{supp}(w)$ contains no consecutive values.

    If $\op{supp}(w) \subseteq [1,n-1]$, then $w$ is free if and only if it is free as an element of $S_n$, and so (by Proposition~\ref{prop:free patterns in type a}) if and only if it classically avoids the patterns $\{231,312,321\}$. 
    On the other hand, if $0 \in \op{supp}(w)$, then $\op{supp}(w) \setminus \{0\}$ is a subset of $[2,n-1]$ that contains no consecutive values. In other words, we have $w(1) = \ol{1}$, and $w(2)\cdots w(n)$ is a string a positive integers classically avoiding the patterns $231$, $312$, and $321$. In particular, for any free $w \in \sb_n$, whatever its support, it suffices to avoid $\{231, 312, 321\}$ and to have the only permissible negative value in the positive indices of $w$ being $w(1) = \ol{1}$. In other words, $w$ must also avoid the patterns $\left\{\ol{1}\,\ol{2}, \
    \ol{2}\,\ol{1}, \
    \ol{2}\,1, \
    2\,\ol{1}, \
    1\,\ol{2}\right\}$.
    
    Conversely, any signed permutation $w \in \sb_n$ that is not free has $i, i+1 \in \op{supp}(w)$ for some $i$, from which it can be shown that $w$ has at least one of the eight listed patterns. Specifically, if $\op{supp}(w)$ contains $0$ and $1$, then $w$ either has multiple negative values in its positive indices (yielding one of the patterns $\left\{\ol{1}\,\ol{2}, \
    \ol{2}\,\ol{1}\right\}$), or only one negative value in those indices but $w(1) \neq \ol{1}$ (yielding one of the patterns $\left\{\ol{2}\,1, \
    1\,\ol{2}, \
    2\,\ol{1}\right\}$).  On the other hand, if $\op{supp}(w)$ does not contain $0$ and $1$ but does contain $i, i + 1$ for some $i \ge 1$ then $w$ classically contains at least one of the patterns $\{231,312,321\}$.
\end{proof}

Our next result shows that freeness is a persistent property. 

\begin{theorem}\label{thm:free is cool}
    The set $\GAV_n(\{231,312,321\})$ of signed permutations in $\sb_n$ that globally avoid the patterns $231$, $312$, and $312$ is equal to the set of free signed permutations. In other words, being free is a persistent property.
\end{theorem}

\begin{proof}
    This argument follows the same outline as previous persistence proofs. Suppose, first, that $w$ classically contains one of the eight patterns listed in the statement of Proposition~\ref{prop:characterizing free signed perms}. Looking at the complete notation, as marked below, each of these patterns globally contains a $231$-, $312$-, or $321$-pattern. 
    $$\begin{array}{lr@{\bigmirror}lcr@{\bigmirror}lcr@{\bigmirror}l}
        & \ol{1} \, \ol{3} \, \ol{2} & \note{2} \, \note{3} \, \note{1}
        & \hspace{.25in}
        & \note{1} \, \ol{2} & \note{2} \, \note{\ol{1}}
        & \hspace{.25in}
        & 2 \, \note{1} & \note{\ol{1}} \, \note{\ol{2}}\\

        \\
        & \ol{2} \, \ol{1} \, \ol{3} & \note{3} \, \note{1} \, \note{2}
        & \hspace{.25in}
        & \ol{1} \, \note{2} & \note{\ol{2}} \, \note{1}
        & \hspace{.25in}
        & \note{2} \, \ol{1} & \note{1} \, \note{\ol{2}}\\
        \\
        & \ol{1} \, \ol{2} \, \ol{3} & \note{3} \, \note{2} \, \note{1}
        & \hspace{.25in}
        & 1 \, \note{2} & \note{\ol{2}} \, \note{\ol{1}}\\
    \end{array}$$
    Thus every non-free signed permutation globally contains a $231$-, $312$-, or $321$-pattern.

    Conversely, consider some $w \in \sb_n$ that globally contains $231$, $312$, or $321$. As in the proof of Theorem~\ref{thm:type b boolean global patterns}, we will explicitly construct a classical pattern characterization of the global pattern class of $\{231,312,321\}$, following the same ideas as \cite[Thm.~2.20]{egge2007}.

    Suppose, first, that $w$ globally contains $231$, and that this global pattern occurrence is in positions $p_1 < p_2 < p_3$, meaning that $w(p_3) < w(p_1) < w(p_2)$. We proceed by cases, based on the signs of the positions $p_i$. 
    For brevity, we leave the detailed checking as an exercise to the reader, but we outline the results here. 
    \begin{itemize}
        \item When $p_1 > 0$ or $p_3 < 0$, this $w$ classically contains at least one element of
        $$\left\{ 231, \
        \ol{2}\,\ol{1}, \ 
        2\,\ol{1}, 
        1\,\ol{2}
        \right\}.$$
        \item When $p_1 < 0 < p_2$ or $p_2 < 0 < p_3$, this $w$ classically contains at least one element of 
        $$\left\{
        \ol{2}\,1, \
        \ol{2}\,\ol{1}, \ 
        2\,\ol{1}, \ 
        1\,\ol{2}
        \right\}.$$
    \end{itemize}
    By Lemma~\ref{lem:egge symmetries}(b), since $\rc{231} = 312$, if $w$ globally contains $312$ then $w$ globally contains $231$, and thus it classically contains at least one of the eight patterns. 

    Finally, suppose that $w$ globally contains $321$, and that this global pattern occurrence is in positions $p_1 < p_2 < p_3$, meaning that $w(p_3) < w(p_2) < w(p_1)$. We again proceed by cases.
    \begin{itemize}
        \item When $p_1 > 0$ or $p_3 < 0$, this $w$ classically contains at least one element of
        $$\left\{ 321, \
        \ol{1}\,\ol{2}, \ 
        2\,\ol{1}, 
        1\,\ol{2}
        \right\}.$$
        \item When $p_1 < 0 < p_2$ or $p_2 < 0 < p_3$, this $w$ classically contains at least one element of 
        $$\left\{
        \ol{2}\,1, \
        \ol{1}\,\ol{2}, \ 
        1\,\ol{2}
        \right\}.$$
    \end{itemize}
    Thus, in whatever way $w$ globally contains a $231$-, $312$-, or $321$-pattern, it classically contains at least one of the patterns listed in the statement of Proposition~\ref{prop:characterizing free signed perms}. Hence, by Proposition~\ref{prop:free patterns in type a}, being free is a persistent property.
\end{proof}

\subsection{Separable elements}

We now consider the class of \emph{separable} (signed) permutations, which has been characterized in terms of both classical and BP-pattern avoidance.  We find that the global description is strictly more attractive than the BP-characterization.

A permutation is separable if it can be obtained from the trivial permutation $1 \in S_1$ by iterating the operations direct-sum and skew-sum.  This condition can be characterized by pattern avoidance: the separable permutations are exactly those that avoid $2413$ and $3142$.  In \cite{gaetz gao}, the notion of separability was extended to all Weyl groups, where it was characterized in terms of BP-avoidance.  For signed permutations, the result is as follows.
\begin{theorem}[part of {\cite[Thm.~5.3]{gaetz gao}}]
A signed permutation $w \in \sb_n$ is separable if and only if $w$ avoids the type-B BP-patterns $2 \, \ol{1}$ and $\ol{2} \, 1$ and the type-A BP-patterns $2413$ and $3142$.
\end{theorem}

As suggested above, the global pattern characterization of these elements is notably simpler. In fact, separability is another persistent pattern property.

\begin{theorem}
A signed permutation $w \in \sb_n$ is separable if and only if $w$ globally avoids the patterns $2413$ and $3142$.
\end{theorem}

This result can be proved in a manner analogous to the results of previous sections, and we omit its details here.

\subsection{Layered permutations}

A permutation is \emph{layered} if it is a direct sum of decreasing permutations, and layered permutations are exactly those that avoid $231$ and $312$. Equivalently, a permutation is layered if and only if it is the longest element of a standard parabolic subgroup of $S_n$.  This latter characterization has an obvious analogue in $\sb_n$: we say a signed permutation is \emph{layered} if and only if it is the longest element of a standard parabolic subgroup of $\sb_n$. For example, the signed permutation $\ol{1}\,\ol{2}\,\ol{3}\,5\,4\,6\,9\,8\,7$ is a layered signed permutation, corresponding to the standard parabolic subgroup generated by $\{\sb_0, \sb_1, \sb_2, \sb_4, \sb_7, \sb_8\}$. 

This is another persistent pattern property.

\begin{theorem}
    A signed permutation $w \in \sb_n$ is layered if and only if $w$ globally avoids the patterns $231$ and $312$.
\end{theorem}

As with separability, we omit the proof of this result, which is similar to the previous proofs in this section.

\subsection{Some properties that are not persistent}

Given the results above, one might hope that type-B characterizations are \emph{always} a globalization of type-A results. This is not, however, the case, as we quickly see with Grassmannian and bigrassmannian elements, which are characterized by having a unique descent and by both the permutation and its inverse having a unique descent, respectively \cite[P0062, P0063, P0064, P0065]{dppa}. Indeed, unsigned permutations of either type must classically avoid $321$ (among other patterns), but $1\, \ol{2}$ is both Grassmannian and bigrassmannian and yet it globally contains a $321$-pattern. Nevertheless, it does seem that Grassmannian and bigrassmannian signed permutations can be characterized by global pattern avoidance.
\begin{conjecture}\label{conj:grassmannian}
    A signed permutation is Grassmannian if and only if it globally avoids the patterns in
    \[
    P := \{4321, 32154, 42153, 43152, 52143, 53142, 
    214365, 315264, 314265, 415263
     \},
    \]
    and a signed permutation is bigrassmannian if and only if it globally avoids the patterns in
    \[
    P \cup \{ p^{-1} : p \in P\}.
    \]
\end{conjecture}
(If Conjecture~\ref{conj:grassmannian} is correct, then the set $P$ is not uniquely determined; for example, each element $p \in P$ listed in Conjecture~\ref{conj:grassmannian} could be independently replaced by $\rc{p}$.)

Similarly, neither rank-symmetry nor self-duality of Bruhat intervals, both of which can be characterized by pattern avoidance in type A \cite[P0005 and P0057]{dppa}, are persistent to type B. Indeed, both properties are thwarted by $\ol{2}\, \ol{1}$ because it has both properties but it globally contains $3412$, one of the forbidden patterns in the type-A characterizations.  It remains to determine whether rank-symmetry or self-duality of Bruhat intervals can be characterized by global avoidance of some other set of patterns.

\section{Enumerative results}\label{sec:enumerations}

It is a foundational result of that field that for any $p \in S_3$, the $p$-avoiding permutations in $S_n$  are enumerated by the $n$th Catalan number $C_n = \binom{2n}{n}/(n+1)$.  The enumerative study of pattern classes has been extended (though not to the same degree) to classical pattern avoidance in signed permutations, as in \cite{simion, mansour-west}.  Enumeration of global pattern classes (under other names) has also been the subject of some prior work, as discussed in Section~\ref{subsec:enum}.  

In Section~\ref{subsec:321}, we offer a new proof of the enumeration of $|\GAV_n(321)|$ that is somewhat more direct and combinatorial than previous proofs.  We then enumerate $\GAV_n(\{p, q\})$ for $p \in S_3$ and $q$ a monotone pattern.  We conclude the section with global pattern results in the vein of the Erd\H{o}s--Szekeres theorem.

\subsection{New perspective on an old result}\label{subsec:321}

As discussed in Section~\ref{sec:fully commutative}, the signed permutations in $\GAV_n(321)$ are precisely Stembridge's top fully commutative elements, and \cite[Prop.~5.9(b)]{stembridge-fc} establishes (using the combinatorics of reduced words, plane partitions, and symmetric functions) that the number of these is the central binomial coefficient $\binom{2n}{n}$ (the \emph{type-B Catalan number} \cite[Rem.~2]{Reiner}).  We are aware of three other proofs of this enumeration: first, in \cite[Thm.~2.17]{egge2007}, Egge found a functional equation satisfied by a generating function for a refined enumeration, taking into account a recursive construction of signed permutations of all sizes globally avoiding $123$, then solved this using the kernel method.  Egge gave a second proof \cite[Cor.~6.3]{egge2010} using direct analogues of the \emph{Robinson--Schensted correspondence} and \emph{Greene's theorem} (as in \cite[Ch.~7]{Stanley 1999}) in the setting of signed permutations (described in more detail below) to express $|\GAV_n(12\cdots k)|$ in terms of \emph{domino tableaux}, then used further bijections to re-express the enumeration in terms of the numbers of (unsigned) permutations avoiding shorter monotone patterns.  Finally, Adenbaum and Elizalde \cite{adenbaum elizalde} showed that the centrally symmetric $321$-avoiding permutations of length $2n$ are in bijection with Dyck paths of semilength $2n$ that are fixed under reversal, whose enumeration is known to be $\binom{2n}{n}$. Their argument views these objects as antichains in the type-B root poset, embedded as the symmetric antichains in the corresponding type-A root poset, by taking such a permutation to the Dyck path whose peaks are given by the permutation's excedances.  

We now offer a new proof of this enumeration, which (like the Adenbaum--Elizalde proof) is bijective and makes use of the \emph{Garfinkle--Barbasch--Vogan (GBV) correspondence} alluded to in the previous paragraph.  In this setting, the standard Young tableaux of Robinson--Schensted are replaced by \emph{(standard) domino tableaux}, which are fillings of a Young diagram $\lambda$ of size $|\lambda| = 2n$ by $n$ dominoes, labelled $1, \ldots, n$, so that the labels increase along each row and column.  The GBV correspondence is a bijection between $\sb_n$ and pairs $(P, Q)$ of domino tableaux of the same shape of size $2n$ (i.e., with $n$ dominoes).\footnote{The GBV correspondence is not the only domino analogue of RS for signed permutations; for example, another bijection was given by Shimozono and White \cite{shimozono white}.  While the Shimozono--White domino algorithm does satisfy \emph{a} version of Schensted's theorem on longest increasing subsequences, it is not well suited to studying global pattern avoidance.}  Barbasch--Vogan \cite{barbasch vogan} and Garfinkle \cite{garfinkle} were interested in this bijection for its applications to Lie theory; a thorough treatment from the combinatorial point of view, including a realization of the GBV correspondence in terms of RS, was given by van Leeuwen \cite[\S4.2]{van leeuwen}, and was implicitly rediscovered by Egge \cite{egge2010}. 
A particular consequence of this connection is the following result (an analogue of Schensted's theorem \cite[Cor.~7.23.11]{Stanley 1999}): a signed permutation $w$ globally avoids $12\cdots k(k + 1)$ if and only if its image $(P, Q)$ under the GBV correspondence consists of a pair of domino tableaux with at most $k$ boxes in the first row; and similarly that $w$ globally avoids $(j + 1)j \cdots 21$ if and only if $P$ and $Q$ have at most $j$ boxes in the first column.  

With this background in hand, we are prepared for our new proof of an old result.  An interactive implementation of the GBV correspondence is available at \cite{devra}. The bijection described in the proof of Theorem~\ref{thm:enumerating GAV_n(123)} is illustrated in Figure~\ref{fig:tableaux}.

\begin{theorem}\label{thm:enumerating GAV_n(123)}
    The number of signed permutations in $\sb_n$ that globally avoid $321$ is $\binom{2n}{n}$. 
\end{theorem}
Of course by symmetry the same result holds for the pattern $123$.

\begin{proof}
Let $B(2n, 2n)$ denote the set of domino tableaux of shape $(2n,2n)$ (i.e., rectangles of two rows with $2n$ boxes in each row).  It is known that $|B(2n, 2n)| = \binom{2n}{n}$; for example, it is not hard to show that associating to each domino tableau of this shape the Dyck path of semilength $2n$ that is fixed under reversal whose $i$th step (for $1 \leq i \leq 2n$) is up if the $i$th domino includes a box in the first row (and is down otherwise) is a bijection.  We give a bijection between $\GAV_n(321)$ and the $B(2n, 2n)$.  Given a signed permutation $w \in \GAV_n(321)$, let $(P,Q)$ be the pair of domino tableaux to which $w$ corresponds under GBV, both of shape $(2n - k, k)$ for some $0 \leq k \leq n$.  Construct a new domino tableau $T$, as follows: rotate $Q$ by $180^\circ$; for each domino in the rotated tableau, replace its label $i$ with $2n + 1 - i$; glue the resulting tableau to $P$.  To reverse the bijection, separate the portion of $T \in B(2n,2n)$ made up of dominoes labelled $1, \ldots, n$ from the portion made up of dominoes labelled $n + 1, \ldots, 2n$; rotate the latter by $180^\circ$ and subtract its labels from $2n + 1$; finally, apply the inverse of GBV to the resulting pair. 
\end{proof}
Alternatively, one could begin by showing that there are $\binom{n}{\lfloor k/2 \rfloor}$ domino tableaux of shape $(2n-k,k)$, and then use the GBV correspondence together with binomial coefficient identities to get the result.

\begin{figure}
\begin{multline*}
w = \ol{3} \, 1 \, 2 
\quad \overset{\text{GBV}}{\longleftrightarrow} \quad
(P, Q) = \left( {\begin{young}[15.5pt][c]
]= 1 & =] 1 & ]= 2 & =] 2 \\
]= 3 & =] 3
\end{young}} \, , \;
{\begin{young}[15.5pt][c]
1 \ynobottom & 2 \ynobottom  & ]= 3 & =] 3 \\
]=] 1 \ynotop & ]=] 2 \ynotop
\end{young}}
\right) \\
 \longleftrightarrow \quad \left( {\begin{young}[15.5pt][c]
]= 1 & =] 1 & ]= 2 & =] 2 \\
]= 3 & =] 3
\end{young}} \, , \;
{\begin{young}[15.5pt][c]
, & , & 5 \ynobottom & 6 \ynobottom \\
]= 4 & =] 4 & ]=] 5 \ynotop & ]=] 6 \ynotop
\end{young}}
\right)
\quad \longleftrightarrow \quad
{\begin{young}[15.5pt][c]
]= 1 & =] 1 & ]= 2 & =] 2 & ]=] 5 \ynobottom & 6 \ynobottom \\
]= 3 & =] 3 & ]= 4 & =] 4 & ]=] 5 \ynotop & ]=] 6 \ynotop
\end{young}}
\end{multline*}
    \caption{The bijection in the proof of Theorem~\ref{thm:enumerating GAV_n(123)}.}
    \label{fig:tableaux}
\end{figure}

\subsection{A pattern of size \texorpdfstring{$3$}{3} and a monotone pattern}
\label{subsec: GAVn(132,123...(k-1))}

In this section, we give formul\ae\ for $|\GAV_n(\{p, q\})|$ when $p \in S_3$ and $q$ is a monotone pattern.  The cases when $p$ is also monotone are either trivial (in the sense that $|\GAV_n(\{p, q\})| = 0$ for sufficiently large $n$) or reduce to the case considered in Section~\ref{subsec:321}.  Thus, up to symmetry, there are only two cases to consider: $\{132, 12\cdots k(k+1)\}$ and $\{132, (k + 1)k\cdots 21\}$.  Our analysis in both cases will be facilitated by an analysis of the structure of $132$-avoiding signed permutations, which was carried out in \cite[Lem.~2.8]{egge2007} and which we repeat here for convenience.

By \cite[Lem.~2.6]{egge2007}, a signed permutation globally avoids $132$ if and only if it avoids both $132$ and $213$, and by Proposition~\ref{prop:bijection between global avoidance and classical avoidance in reverse-complement invariant permutations} such signed permutations correspond (under $\embed$) to the rc-invariant $132$- and $213$-avoiding permutations in $S_{2n}$.  It is well known that the $132$- and $213$-avoiding permutations in $S_n$ are the \emph{colayered permutations}, namely the permutations $w$ for which there are indices $0 = i_0 < i_1 < \dots < i_k = n$ such that 
\begin{multline}\label{eq:colayered}
n = w_{i_1} > w_{i_1 - 1} > \dots > w_1 > {} 
w_{i_2} > w_{i_2 - 1} > \dots > w_{i_1 + 1} >  {}\\ 
w_{i_3} > w_{i_3 - 1} > \dots > w_{i_2 + 1} >  \dots > {}\\ 
 w_{i_{k - 1}} > w_{i_{k - 1} - 1} > \dots > w_{i_{k - 2} + 1} > 
w_n > w_{n - 1} > \dots > w_{i_{k - 1} + 1} = 1.
\end{multline}
Such permutations are in natural correspondence with compositions of $n$, where the permutation satisfying \eqref{eq:colayered} corresponds to the composition $(i_1 - i_0, i_2 - i_1, \ldots, i_k - i_{k - 1})$ that records the lengths of the increasing runs of $w$.  

It follows immediately from the preceding considerations that the set $\GAV_n(132)$ of globally $132$-avoiding signed permutations in $\sb_n$ is in natural bijection with the set of \emph{palindromic} compositions of $2n$, i.e., the compositions $(a_1, \ldots, a_k)$ such that $a_i = a_{k + 1 - i}$ for $i = 1, \ldots, k$. There are $2\cdot 2^{n-1} = 2^n$ of these: take any of the $2^{n-1}$ compositions of $n$, mirror it and then choose whether or not to consolidate the two central parts.

\begin{example}\label{ex:colayered}
    The signed permutation $1\, \ol{4} \, \ol{3} \, \ol{2} \, \ol{6} \, \ol{5}$ belongs to $\GAV(132)$.  Its image under the inclusion map $\embed$ is $11\,12\,8\,9\,10\,6\,7\,3\,4\,5\,1\,2 \in S_{12}$, graphed in Figure~\ref{fig:colayered}, which is colayered and corresponds to the palindromic composition $(2, 3, 2, 3, 2)$.
    \begin{figure}[htbp]
        \begin{tikzpicture}[scale=.25]
            \foreach \x in {1,2,3,4,5,6,7,8,9,10,11,12} {
                \draw (\x,1) -- (\x,12);
                \draw (1,\x) -- (12,\x);
            }
            \foreach \x in {(1,11),(2,12),(3,8),(4,9),(5,10),(6,6),(7,7),(8,3),(9,4),(10,5),(11,1),(12,2)}{
                \fill \x circle (7pt);
            }
        \end{tikzpicture}
        \caption{The colayered permutation $11\,12\,8\,9\,10\,6\,7\,3\,4\,5\,1\,2 \in S_{12}$.}\label{fig:colayered} 
    \end{figure}
\end{example}

We now consider the addition of a further constraint on the longest increasing or decreasing subsequence of $w$. 
In the case of an increasing forbidden subsequence, we show that the global avoidance classes share initial values and a recurrence relation with the Fibonacci-like sequence 
    \[
        f^{(k)}_i := \begin{cases}
        \displaystyle \sum_{j = 1}^k f^{(k)}_{i - j} & \text{if } i > k,\\
        1 & \text{if } i = \lfloor k/2 \rfloor + 1, \text{ and} \\
        0 & \text{otherwise},
        \end{cases}
    \]
defined for positive integers $i$ and $k$. This sequence has a number of other combinatorial interpretations; see, for example, \cite{prodinger}.  We begin with a technical lemma.

\begin{lemma}\label{lem:fibonacci-like sequence}
    For positive integers $k$, 
    $$f^{(k)}_{n+k+1} = \begin{cases}
        2^n & \text{if } 0 \le n \le \lfloor k/2 \rfloor, \text{ and}\\
        2^n - 2^{n-\lfloor k/2\rfloor - 1} & \text{if } \lfloor k/2 \rfloor < n < k.
    \end{cases}$$
\end{lemma}

\begin{proof}
    We proceed by induction on $n$, starting with the base case $f^{(k)}_{k + 1} = 0 + \cdots + 0 + 1 + 0 + \cdots + 0 = 2^0$.
    If $1 \le n \le \lfloor k/2 \rfloor$, then 
    $$
        f^{(k)}_{n + k + 1}  = \sum_{j = 1}^{n} f^{(k)}_{n + k + 1 - j} + \sum_{j=n + 1}^{k} f^{(k)}_{n + k + 1 - j}
        = \sum_{t = 0}^{n - 1} f^{(k)}_{t + k + 1} + \sum_{i = n + 1}^{k} f^{(k)}_{i}
        = \left(\sum_{t=0}^{n-1} 2^{t}\right) + 1 
        = 2^n. 
    $$
    For $\lfloor k/2 \rfloor < n < k$, on the other hand, we have
   \begin{align*} 
        f^{(k)}_{n + k + 1} & \ = \ \sum_{j = 1}^{n}  f^{(k)}_{n+k+1-j} +  \sum_{j = n + 1}^{k}  f^{(k)}_{n+k+1-j} 
        \ = \ \sum_{t = 0}^{n - 1}  f^{(k)}_{t+k+1} +  \sum_{i = n + 1}^{k}  f^{(k)}_{i} = \sum_{t = 0}^{n - 1}  f^{(k)}_{t+k+1}\\
        & \ = \  \sum_{t = 0}^{\lfloor k/2\rfloor}  f^{(k)}_{t+k+1} \hspace{.05in} + \!\!\! \sum_{t = \lfloor k/2\rfloor + 1}^{n - 1}  \!\!\! f^{(k)}_{t + k+1} 
        \ = \ \sum_{t = 0}^{\lfloor k/2\rfloor}  2^t  \hspace{.05in} + \!\!\! \sum_{t = \lfloor k/2\rfloor + 1}^{n - 1} \!\!\! \left(2^t - 2^{t - \lfloor k/2\rfloor - 1}\right)
          \\
        & \ = \  \sum_{t = 0}^{n - 1} 2^t \hspace{.05in} -\!\!\! \sum_{t = \lfloor k/2\rfloor + 1}^{n - 1} \!\!\! 2^{t - \lfloor k/2\rfloor - 1} 
        \ = \ (2^n - 1) - \left(2^{n - \lfloor k/2\rfloor - 1} - 1\right) 
        \ = \ 2^n - 2^{n - \lfloor k/2 \rfloor - 1}. && 
    \end{align*}
This completes the proof.
\end{proof}

The next two results enumerate global avoidance of $132$ and a monotone pattern.

\begin{theorem}
\label{GAVn(132,123...(k-1))}
    For $k \geq 1$, the number of signed permutations in $\sb_n$ that globally avoid the patterns $132$ and $12 \cdots k(k + 1)$ is 
$f^{(k)}_{n + k + 1}$. 
\end{theorem}

\begin{proof}
We show that enumeration of the signed permutations in question agrees with the initial values and recurrence relation of the sequence discussed in Lemma~\ref{lem:fibonacci-like sequence}.

Obviously $\GAV_n(\{132, 12\cdots (k + 1)\}) = \GAV_n(132)$ if $2n < k + 1$; since $|\GAV_n(132)| = 2^n$ for all $n \geq 0$ (taking by convention $|\sb_0| = 1$), this gives the required initial values for $n \leq \lfloor k/2\rfloor$.  We now consider larger values of $n$.
Suppose that, under the correspondence described at the beginning of this section, a signed permutation $w$ in $\GAV_n(132)$ corresponds to the palindromic composition $(a_1, \ldots, a_m)$.  We claim that the length of the longest increasing subsequence in $w$ is $\max\{a_1, \ldots, a_m\}$, and consequently that $w$ globally avoids $12\cdots k(k + 1)$ if and only if $a_i \leq k$ for all $i$.  Indeed, this holds because (by the structure of colayered permutations) every pair of entries of $w$ that belong to different increasing runs occur in decreasing order in $w$, and so any increasing subsequence must belong entirely to one of the increasing runs.
It follows immediately that $|\GAV_n(\{132, 12\cdots (k + 1)\})| = 2^n - 2^{n - \lfloor k/2\rfloor - 1}$ for $k/2 < n < k$, since the palindromic compositions of $2n$ in this range that contain a part of size $k + 1$ or larger have the form $(a_1, \ldots, a_\ell, z, a_\ell, \ldots, a_1)$, where $(a_1, \ldots, a_\ell)$ is any composition of some integer less than or equal to $n - \lfloor k/2\rfloor - 1$ (including possibly the empty composition of $0$), and there are precisely $2^{n - \lfloor k/2\rfloor - 1}$ of these.

Now suppose that $n > k$.  Then any palindromic composition $(a_1, \ldots, a_m)$ of $2n$ in which all parts are at most $k$ has $m > 2$, and the trimmed composition $(a_2, \ldots, a_{m - 1})$ is a palindromic composition of $2n - 2a_1$ with all parts at most $k$.  By induction, the number of these is $f^{(k)}_{n - a_1 + k + 1}$, and $a_1$ may be any element of $\{1, 2, \ldots, k\}$.  Thus 
\[
\left|\GAV_n(\{132, 12\cdots (k + 1)\})\right| = \sum_{j = 1}^k f^{(k)}_{n - j + k + 1} = f^{(k)}_{n + k + 1}.  \qedhere
\] 
\end{proof}

Substituting $k=2$ into Theorem~\ref{GAVn(132,123...(k-1))} recovers Egge's enumeration of $|\GAV_n(\{123,132\})|$ by the standard Fibonacci numbers. The counterpoint to Theorem~\ref{GAVn(132,123...(k-1))}, pairing the other monotone pattern with $132$, yields the following enumeration.

\begin{theorem}
  For $k \geq 1$, the number of signed permutations in $\sb_n$ that globally avoid the patterns $132$ and $(k + 1)k \cdots 21$ is
  \[
  \sum_{j = 1}^{k} \binom{n - 1}{\lfloor (j - 1)/2\rfloor}.
  \]
\end{theorem}

\begin{proof}
Suppose that, via the correspondence described earlier in this section, the signed permutation $w \in \GAV_n(132)$ corresponds to the palindromic composition $(a_1, \ldots, a_m)$.  We claim that the length of the longest decreasing subsequence in $w$ is $m$. Indeed, each decreasing subsequence can contain at most one entry from each of the $m$ increasing runs, and any such selection does give a decreasing subsequence in $w$ of length $m$ (by the colayered property).  It follows immediately that $|\GAV_n(\{132, (k + 1)k\cdots 21\})|$ is equal to the number of palindromic compositions of $2n$ having at most $k$ parts. 

For even $j$, there is a bijection between palindromic compositions of $2n$ into $j$ parts and compositions of $n$ into $j/2$ parts: $(a_1, \ldots, a_{j/2}, a_{j/2}, \ldots, a_1) \leftrightarrow (a_1, \ldots, a_{j/2})$. 
For odd $j$, there is a bijection between palindromic compositions of $2n$ into $j$ parts and compositions of $n$ into $(j + 1)/2$ parts, with
$(a_1, \ldots, a_{(j - 1)/2}, a_{(j + 1)/2}, a_{(j - 1)/2}, \ldots, a_1)$
(whose middle term must be even) corresponding to  $(a_1, \ldots, a_{(j - 1)/2}, \frac{1}{2}a_{(j + 1)/2})$.  There are $\binom{n - 1}{(j + 1)/2 - 1} = \binom{n - 1}{\lfloor (j - 1)/2\rfloor}$ compositions in each case.

The result now follows by summing over the valid possibilities for $j$.
\end{proof}

\subsection{An Erd\H{o}s--Szekeres-type result}\label{sec:erdos szekeres}

The increasing permutation $1 2\cdots (n - 1) n$ only classically contains other increasing permutations, and the decreasing permutation $n (n - 1) \cdots 2 1$ only classically contains other decreasing permutations. Thus, for a set $P$ of patterns, the collection of permutations in $S_n$ classically avoiding $P$ will be nonempty unless $P$ contains both an increasing and decreasing pattern.  The Erd\H{o}s--Szekeres theorem provides a formula for the largest size of an unsigned permutation that classically avoids the patterns $\{12\cdots k(k + 1), (j+1)j\cdots21\}$, and the number of permutations of this maximum size can be computed easily using the RS correspondence and Schensted's theorem, as discussed in Section~\ref{subsec:321}.

\begin{proposition}\label{prop: ES1}
If $w \in S_n$ avoids $12\cdots k(k + 1)$ and $(j + 1)j\cdots 21$, then $n \leq kj$.  Moreover, the number of permutations in $S_{kj}$ that avoid these patterns is the square of the number of standard Young tableaux of shape $(\underbrace{k, k, \ldots, k}_{j \text{ copies}})$.
\end{proposition}

It follows from the \emph{hook-length formula} (see, e.g., \cite[\S7.21]{Stanley 1999}) that the number of tableaux referenced in Proposition~\ref{prop: ES1} is the $k$th \emph{$j$-dimensional Catalan number} \cite[A060854]{OEIS}
\[   
C_{(j,k)} = \frac{(kj)! \cdot \left(1! \cdot 2! \cdots (j - 1)!\right) \cdot \left(1! \cdot 2! \cdots (k - 1)!\right) }{1! \cdot 2! \cdots (k + j - 1)!}.
\]
Replacing permutations with signed permutations, and RS with GBV, leads to the following analogue of the previous proposition. 

\begin{proposition}\label{prop: ES2}
If $w \in \sb_n$ globally avoids $12\cdots k(k + 1)$ and $(j + 1)j\cdots 21$, then $n \leq \frac{kj}{2}$.  Moreover, the number of permutations in $\sb_{\lfloor kj/2\rfloor}$ that globally avoid these patterns is 
the square of the number of domino tableaux whose shape is the $j$-tuple 
$$\begin{cases}
    (k,k,\ldots, k,k) & \text{if $kj$ is even, and}\\
    (k,k,\ldots, k,k-1) & \text{if $kj$ is odd.}
\end{cases}$$
\end{proposition}

Moreover, it follows from \cite[Thm.~4.3]{stembridge-dt} that the number of domino tableaux referenced in Proposition~\ref{prop: ES2} is the product of a binomial coefficient with the numbers of standard Young tableaux of two smaller rectangles (which are themselves given by higher-dimensional Catalan numbers).

\section{Open questions}\label{sec:open}

The breadth and depth of the results proved above suggest several directions for further research. In this section, we propose a selection of them. To give additional context to these questions, we group some of them by topic.

\subsection{Characterization and persistence}

Above, we conjectured (Conjecture~\ref{conj:grassmannian}) that Grassmannian and bigrassmannian signed permutations can be characterized by global pattern avoidance, and we raised the question of whether the same is true for rank-symmetry or self-duality of Bruhat intervals.  One can ask the same question about other natural families of signed permutations, such as Lambert's theta-vexillary permutations \cite{lambert} \cite[P0054]{dppa}.

Another natural question along these lines concerns the \emph{$2$-boolean} elements of a Coxeter group, introduced by Gao and H\"anni in \cite{gao hanni boolean} and defined as those elements whose reduced decompositions contain each generator at most twice. Gao and H\"anni show that in the symmetric group this property is characterized by the avoidance of the patterns $3421$, $4312$, $4321$, and $456123$.  Is this property persistent?  

More broadly, is there some characterization of which (signed) permutation properties are persistent?

\subsection{Enumeration}

There are now several proofs that $|\GAV_n(123)| = \binom{2n}{n}$, all of which are somewhat involved or rely on complicated machinery.  Is there a direct argument using the fact that $\binom{2n}{n} = (n + 1) \cdot C_n$ and that $C_n$ is the number of classical $123$-avoiders in $S_n$?  What about Simion--Schmidt-style bijections between $|\GAV_n(123)|$ and the signed permutations considered by Reading in \cite[Thm.~7.8]{reading}?

Eric Egge and Ben Adenbaum have separately explained to us that the separable signed permutations of size $n$ (those in $\GAV_n(\{2413, 3142\})$) are enumerated by the large Schr\"oder numbers \cite[A006318]{OEIS}.  Large Schr\"oder numbers also enumerate separable permutations, with an indexing shift: the number of separable signed permutations in $\sb_n$ is equal to the number of separable permutations in $S_{n + 1}$.  This parallels the situation for smooth elements described in Section~\ref{subsec:enum}. Is this phenomenon, that the enumerations for $\sb_n$ and $S_{n + 1}$ globally avoiding a set of patterns agree, more widespread?

Computational data suggest that the signed permutations in $\GAV_n(2413)$ (equivalently, in $\GAV_n(3142)$) are enumerated by \cite[A115197]{OEIS}, a sequence of convolutions of certain generalized Catalan numbers. Is this correct? If so, why?  What other familiar sequences enumerate global pattern classes?

In a different vein, what kind of Wilf-equivalence occurs in the context of global pattern containment? Initial steps in this direction appear in Section~\ref{sec:global} above, but perhaps more can be said.

\subsection{Other}

For each set $P$ of global patterns, there is an associated set $\gl(P)$ of classical patterns (defined in \cite[Thm.~2.20]{egge2007}) that characterize its avoidance class.  Generally $\gl(P)$ is much larger than $P$.  How much can be said about how large or small $\gl(P)$ can be in terms of the size $|P|$?  In the spirit of some of the characterization ideas listed previously, can we classify the sets of signed patterns $Q$ that admit a set $P$ of global patterns such that $Q=\gl(P)$? 

In another direction, we wonder what can be said about higher-order containment of global patterns. For example, can one characterize the signed permutations that globally contain $k$ distinct copies of $213$? Grouping by the statistic $k$, how does this partition $\sb_n$?

Finally, are there any interesting manifestations of global pattern avoidance in type D (the even-signed permutations), or in the affine Coxeter groups of classical type (i.e., the other \emph{George groups}), beyond those already known \cite{lewis tenner global, Green}?

\section*{Acknowledgements}

We appreciate the comments, suggestions, and pointers from Ben Adenbaum, David Anderson, Sara Billey, Eric Egge, Yibo Gao, Nathan Reading, and anonymous referees.

\end{document}